\documentclass[a4paper, final, 12pt]{article}
\usepackage{eurosym}
\usepackage{amsmath, amssymb, latexsym, amscd, amsthm,amsfonts,amstext}
\usepackage[mathscr]{eucal}
\usepackage{graphicx}
\usepackage{subfig}
\usepackage{float}
\usepackage{color}
\usepackage{hyperref}
\usepackage[utf8]{inputenc}
\usepackage[english]{babel}
\usepackage[pagewise]{lineno}

 \renewcommand{\baselinestretch}{1}

\numberwithin{equation}{section}
\begin{document}
%
\title{Using the Newton-Raphson Method with Automatic Differentiation to Numerically Solve the Implied Volatility of Stock Options via the Binomial Model}
%
%
%


\author{Wanchaloem~Wunkaew$^{\ast }$,
Yuqing~Liu$^{\circ }$ \\
and Kirill V. Golubnichiy $^{\bullet }$ \and $^{\ast }$Department of Applied Mathematics and Statistics, \and Johns Hopkins University,
MD 21218, USA \and $^{\circ }$Department of Mathematics \and University of
Washington, Seattle, WA 98195, USA \and $^{\bullet }$Department of Mathematics and Statistics, \and Texas Tech University, Lubbock, TX 79409, USA \and Emails: wwunkae1@jh.edu,
\and yliu876@uw.edu, \and %
kgolubni@ttu.edu}
\date{}
\maketitle


\begin{abstract}

In the paper by Michael Klibanov et al. \cite{Kirill}, a novel method is proposed to calculate the implied volatility of European stock options as a solution to the ill-posed inverse problem for the Black-Scholes equation. Additionally, the paper introduces a trading strategy based on the difference between the implied volatility of the option and the volatility of the underlying stock. Besides the Black-Scholes equation, the Binomial model is another method used to price European options, and the implied volatility can also be calculated using this model. In this paper, we apply the Newton-Raphson method in conjunction with Automatic Differentiation to numerically approximate the implied volatility of a stock option using the Binomial model. We provide an explanation of the mathematical model, the methodology, and the results from our tests using simulated data from the Geometric Brownian Motion model and the Binomial model itself, as well as data from the U.S. market from 2020 to 2022.

\end{abstract}


\section{Introduction} 

An option is a financial instrument that gives the holder the right to buy or sell an underlying asset at a specific strike price by a specific date. As these contracts are traded on exchanges, various metrics are used to analyze the behavior of options, with implied volatility being one of the key metrics. Implied volatility measures the expected fluctuations in option prices and has been demonstrated to support decision-making in option trading \cite{Kirill}. In this paper, we present a numerical method to approximate implied volatility by solving an option pricing model known as the Binomial model (or Lattice model) in reverse.

An option to buy is called a "call option," while an option to sell is known as a "put option." Options are further categorized based on the period in which they can be exercised. An American option can be exercised at any time up to its maturity, whereas a European option can only be exercised on the maturity date.

The Black-Scholes model, also known as the Black-Scholes-Merton model, is among the most widely used models in modern financial theory for estimating the theoretical value of options. Developed in 1973 by Fischer Black, Robert Merton, and Myron Scholes, the model, first introduced in \cite{BlackScholes}, uses current stock prices, expected dividends, the option’s strike price, expected interest rates, time to expiration, and expected volatility to calculate the theoretical value of an option contract. The model can be expressed as a parabolic partial differential equation, solvable explicitly. However, despite the model's simplicity, solving for implied volatility using this equation is challenging due to its ill-posed nature \cite{Kirill}.

Another widely used model for option valuation is the Binomial Option Pricing Model \cite{COX1979229}. Initially proposed by William Sharpe in 1978 and later formalized by Cox, Ross, and Rubinstein in 1979, this model uses a discrete tree structure to model the underlying variables of an option over time. In this paper, we discuss using the Binomial Option Pricing Model to calculate implied volatility, detailing the methodology and presenting numerical results from our experiments.

In \cite{Kirill}, a new mathematical method is proposed to accurately calculate the implied volatility of a European stock option. The paper also derives a trading strategy based on the difference between the implied volatility of an option and the volatility of its underlying stock. Specifically, let $\sigma$ represent the volatility of a stock, and $\hat{\sigma}$ represent the volatility of the corresponding option. Let $S$ denote the price of the stock, $f(S,K)$ represent the payoff or value function of the option at maturity time $T$, and $u(S,t)$ denote the price of the option with $t$ representing the time until maturity.

According to the assumptions in the paper, the following equation is satisfied \cite{Kirill}:
\begin{align*}
    \frac{\partial u(S,t)}{\partial t} &= \frac{\hat{\sigma} ^2}{2} S^2 \frac{\partial^2 u(S,t)}{\partial S^2},\qquad S > 0\\
     u(S, 0) &= f(S)
\end{align*}

Let $K$ be the strike price. The payoff function is $f(S) = \max(S - K, 0)$, and $u(S,t)$ is given by the Black-Scholes formula \cite{BlackScholes}:

\begin{align*}
    u(S,t) = S \Phi (\Theta_+ (S,t)) - e^{-rt} K \Phi (\Theta_{-} (S,t))
\end{align*}

where the risk-free rate $r$ is assumed to be $0$.
The stochastic equation of the Geometric Brownian motion for the stock price $S_t$ with volatility $\sigma$ and zero drift is given by

\begin{align*}
    dS_t = \sigma S_t d W_t
\end{align*}
Since the expectation of the Wiener process $dW$ is 0, substituting the above equations into the Ito formula gives the expected value of the increment of the option price over an infinitesimal time interval a

\begin{align*}
    \frac{(\sigma^2 - \hat{\sigma} ^2)}{2} S^2 \frac{\partial^2 u(S,t)}{\partial^2 S}dt
\end{align*}

Since $\frac{\partial^2 u(S,t)}{\partial S^2}$, denoted as the Greek letter $\Gamma$, is proven to be non-negative, and $S^2$ is greater than 0, the expected direction of the option price depends on $\sigma - \hat{\sigma}$. Specifically, under these assumptions, we expect the option price in the market to increase if $\sigma - \hat{\sigma} > 0$, meaning the volatility of the underlying asset is greater than the implied volatility of the option. After deriving the equation, the paper also examines the validity of the mathematical model and generates a possible winning strategy based on it. Testing the hypothesis yields an accuracy of 55.57\%.

The results from \cite{Kirill} reinforce the role and importance of implied volatility in financial mathematics. However, the Black-Scholes equation, on which the paper relies, is not the only method used to price options, and it is applicable only to European options.

In this paper, we first present the algorithm of the Binomial model and explain why it can be used to find the implied volatility of an option's price by demonstrating the convergence from the Black-Scholes equation to the Binomial model. We then discuss the method of calculating implied volatility inversely via the Binomial model, explaining the necessary conditions and parameters for the model.

Next, we test our model by using the Binomial model to calculate the implied volatility of an option price generated from the Geometric Brownian motion model with $\mu = 0$ and $\sigma = 0.2$. We compare the results from the Binomial model to those from the Black-Scholes model and analyze the relative difference in the accuracy of the two models.

\section{Mathematical Model}
\subsection{Binomial Model}
The Binomial Model is based on the key assumption that the price of the underlying asset, in this case the stock, can move in only two possible directions: either up or down in the next time step. These stock prices can be estimated using Geometric Brownian motion:
\begin{align*}
    d S_t &= \mu S_t dt + \sigma S_t d W_t
\end{align*}

where $S_t$ is the stock price at time $t$, $\mu$ is the drift, $\sigma$ is the volatility of the stock, and $W_t$ is a Wiener process.
Note that in the next time step $\delta t$, the Geometric Brownian motion can be approximated by two possible values: \cite{Ross}:
\begin{align*}
    S_{t + \delta t} &= S_t e^{\sigma \sqrt{\delta t}} \quad \text{ or } \quad S_{t + \delta t}= S_t e^{-\sigma \sqrt{\delta t}}
\end{align*}

Assume that the next time step, at time $T$, is the maturity of a call option with strike price $K$. The value of the option at that step can be determined by the following function:
\begin{align*}
    V(T) &= \max(0, S(T) - K) 
\end{align*}

where $S(T)$ is the price of the stock at time $T$, the maturity of the call option.

Since we assume that the stock price can move in only two directions in the next step, we can calculate the value of the call option for each scenario as follows:
\begin{align*}
    C^+ &= \max(0, S(0) e^{\sigma \sqrt{\delta t}} - K) =  \max(0, S_0 e^{\sigma \sqrt{T}} - K)\\
    C^- &= \max(0, S(0) e^{-\sigma \sqrt{\delta t}} - K) = \max(0, S_0 e^{-\sigma \sqrt{T}} - K) \\
\end{align*}
Then, A single-step Binomial Model is formulated as shown below:
\begin{align*}
    C = e^{-r{\delta}t}\left((\frac{1}{2}+ \frac{r \sqrt{\delta t}}{2 \sigma}) C^{+} + (\frac{1}{2} - \frac{r \sqrt{\delta t}}{2 \sigma}) C^{-}\right)
\end{align*}

where $C$ is the value of the option,
$r$ is the risk-free rate,
$t$ is the time step (in years),
$\sigma$ is the volatility of the underlying asset.

This model can be extended into multiple steps, which is referred to as the Lattice Model. In this case, the option values at any possible branch of the tree can be determined using the value functions mentioned above. However, in the intermediate steps, the values of options at future steps are required to calculate the option value at the current step. This process involves solving the option values step by step, working backwards from the final step to determine the option value at the current time.

\subsection{Numerical Methods}

\subsubsection{Newton-Raphson method}The Newton-Raphson method is widely used to numerically solve for the roots of differentiable functions \cite{newton}. The algorithm is expressed by the following formula: \begin{align*} x_{n+1} &= x_n - \frac{f(x_n)}{f'(x_n)} \end{align*} where $x_n$ for $n \in \mathbb{N}$ is the approximated root at the $n^{\text{th}}$ step, $x_0$ is the initial guess, and $f$ is the function. Under stable conditions, $x_n$ will converge to a root of $f$.
\subsubsection{Automatic Differentiation} As shown in the previous subsection, the Newton-Raphson method requires the derivative of $f$ to calculate a root for $f$. There are several methods to approximate the derivative. One intuitive approach is to directly calculate it using theorems from calculus, a process known as symbolic differentiation. However, due to the complexity introduced by formulas such as the product rule and chain rule, this method can lead to overly complex expressions. While this symbolic approach may be feasible for single-step binomial models, it becomes impractical for multi-step lattice models. An alternative numerical method that is widely used is the finite difference method:

\begin{align*}
    f'(x) &\approx \frac{f(x+\Delta x) - f(x)}{\Delta x}
\end{align*}

For some small $\Delta x > 0$, the finite difference method can be used to approximate the derivative for the Newton-Raphson method, and when modified, this approach is known as the Secant Method.

In this project, we apply a widely-used differentiation method called Automatic Differentiation (AD). AD calculates the derivative of a function immediately after each component of the function is computed. The derivative of each component is then stored along with the result of that computation \cite{AD}.

Automatic Differentiation is extensively used in machine learning algorithms, such as backpropagation. Given the structure of the Binomial/Lattice Model, where the option price is calculated through multiple steps in a tree, we propose using Automatic Differentiation to calculate the derivative of the model function.

\section{Methodology} 
In this project, we approximate the implied volatility of a European stock option using the Binomial/Lattice model.
\subsection{Calculating implied volatility inversely via binomial model}
As mentioned in the previous section, to determine the current value of an option, parameters such as the current stock price and risk-free rate are required. These parameters can be obtained from the exchange market, where the option price is also determined.

Now, assuming the known parameters are constants, we can write the formula as follows:
\begin{align*}
    C &= f(\sigma)
\end{align*}
The $\sigma$ that satisfies the equation above is the implied volatility of the corresponding stock option. It is important to note that in this paper, a 10-layer lattice model is applied. The number of layers is determined by computational limitations. The risk-free rate is based on the 10-year treasury yield for real market data, and we assume the rate to be $0$ for the generated dataset. Other parameters are obtained from the real market. The option price $C$ is known in the exchange market, so $\sigma$ can be determined by finding the root of the equation:
\begin{align*}
   f(\sigma) - C &= 0
\end{align*}
We denote the root of the equation above by $\hat{\sigma}$.
To find the root of the equation, we apply the Newton-Raphson method with an initial value of $x_0 = 0.2$ and a tolerance rate of $10^{-5}$. The iteration stops at step $n$ when $|f(x_n) - f(x_{n-1})| < 10^{-5}$. The maximum number of iterations is set to 100. The derivative of $f$ is calculated using Automatic Differentiation at each iteration.
\subsection{Testing} To assess the efficiency of the model, we compare the performance of the proposed method, where $f'(x)$ is computed via automatic differentiation, with the same method where $f'(x)$ is computed via the secant method. For the secant method, we set $\Delta x$ to $10^{-5}$. We study two performance metrics: the number of root convergences and the rate of convergence. Furthermore, we examine these metrics on different datasets as mentioned below.

\subsubsection{Testing the Numerical Method} Since the model used in this project involves multiple step calculations, it is challenging for the numerical methods mentioned above to accurately approximate the implied volatility. Therefore, to assess the efficiency of the model, we tested the algorithms with simulated data.

Similar to \cite{Kirill}, the underlying stock is simulated using a Geometric Brownian motion model with $\mu = 0$ and $\sigma = 0.2$. The prices of options are calculated using the same binomial model, with 90 days to maturity and varying volatilities. After obtaining the simulated option prices, we inversely solve for the implied volatility using the methods mentioned above.

\begin{figure}[h]
\par
\begin{center}
\includegraphics[width =1.00\textwidth]{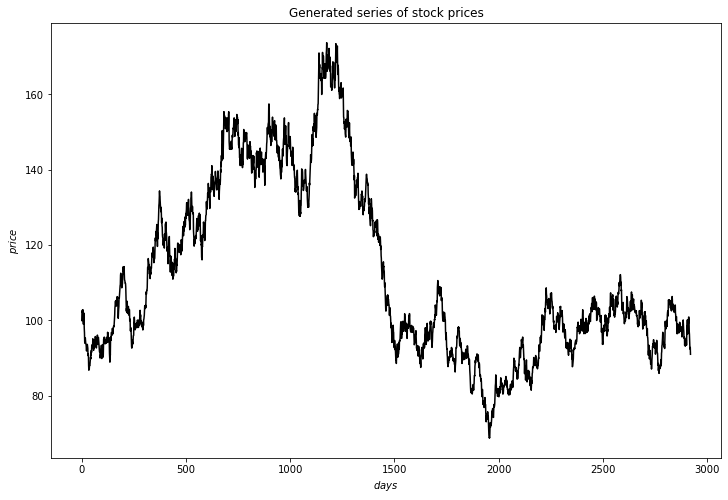}
\end{center}
\par
\caption{Simulated stock price based on Geometric Brownian motion.}
\end{figure}

\begin{figure}[h]
\par
\begin{center}
\includegraphics[width =1.00\textwidth]{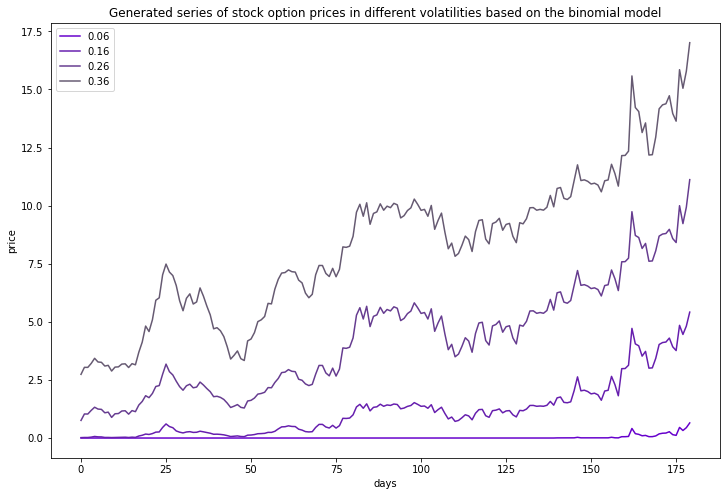}
\end{center}
\par
\caption{An example of simulated option prices generated from the stock prices above and the Binomial model with four different volatilities, a strike price of 150, and a fixed maturity of 180 days.}
\end{figure}

\subsubsection{Testing with Real-World Data} After testing the method with synthetic data, we assess the efficiency of the model using real-world data. In this project, we apply the algorithm to 1917 European options traded in the market from 2020 to 2022. The last prices of each option and its underlying asset are used in this evaluation.
\section{Numerical Findings}
The tables below show the number of data points where the implied volatility converges. In addition, we compute the average number of iterations required for the root to converge.

\begin{center}
\captionof{table}{Generated Data (of total 1350 data)}
\begin{tabular}{|c|c|c|}
    \hline
    & Secant & Autodiff \\ 
    \hline
    Number of Convergences & 1322 & 1308 \\ 
    \hline
    Average Number of Iterations & 3.670 & 3.352 \\ 
    \hline
\end{tabular}

\end{center}

\begin{center}
\captionof{table}{Real Data (of total 1917 data)}
\begin{tabular}{|c|c|c|}
    \hline
    & Secant & Autodiff \\ 
    \hline
    Number of Convergences & 1431 & 1425 \\ 
    \hline
    Average Number of Iterations & 6.304 & 4.700 \\ 
    \hline
\end{tabular}

\end{center}

\subsection{Generated Data}

\clearpage

\begin{figure}[tph]
\centering
\begin{tabular}{cc}
{\includegraphics[width = 0.5\textwidth, height =
0.55\textwidth]{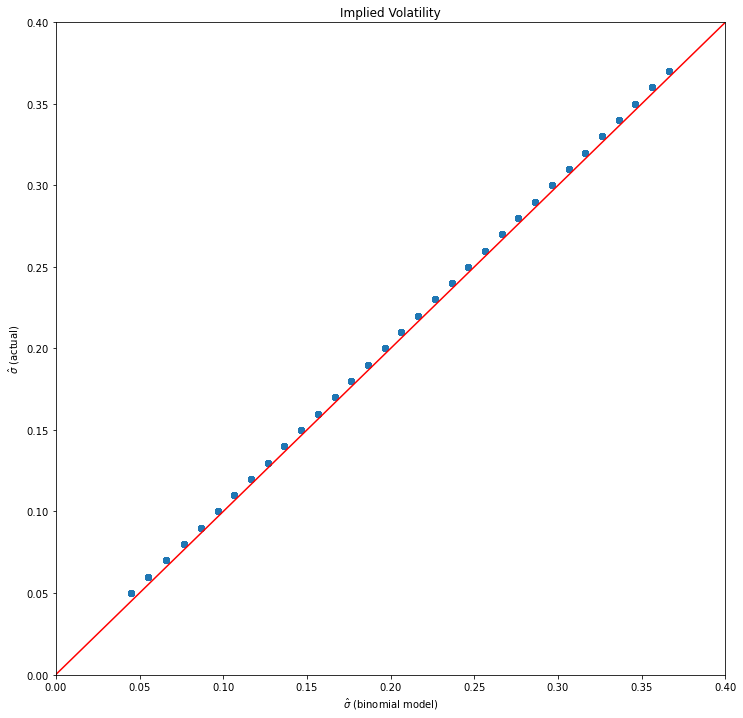}} & {%
\includegraphics[width = 0.5\textwidth,
height = 0.55\textwidth]{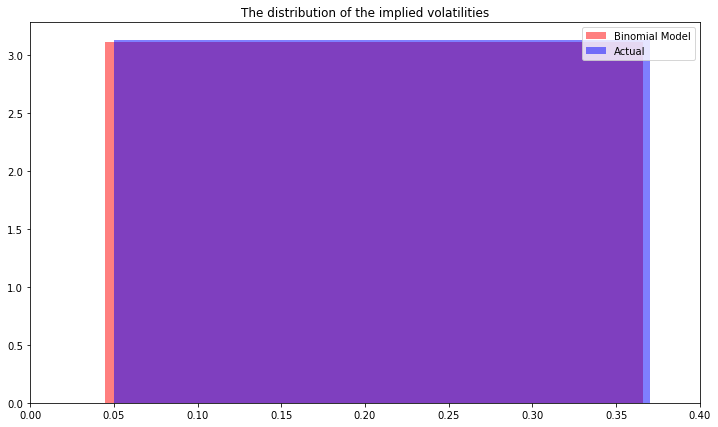}} \\ 
\end{tabular}
\caption{ (left) The scatter plot of implied volatilities calculated from the binomial
model and their corresponding actual values from the dataset (Generated Data); (right)
The implied volatilities distributions (Generated Data). }
\end{figure}

\subsection{Stock options selected from American markets from 2017-
2018}

\clearpage

\begin{figure}[tph]
\centering
\begin{tabular}{cc}
{\includegraphics[width = 0.5\textwidth, height =
0.55\textwidth]{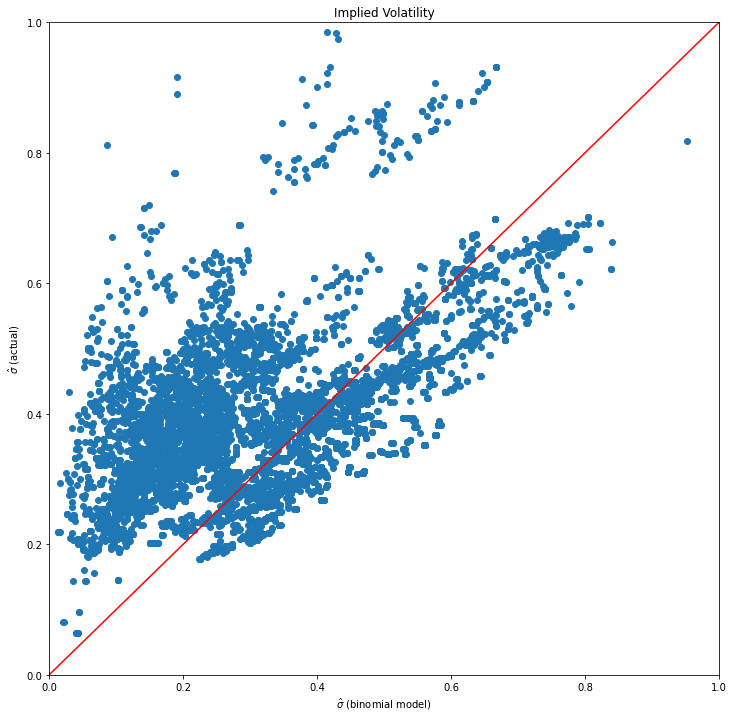}} & {%
\includegraphics[width = 0.5\textwidth,
height = 0.55\textwidth]{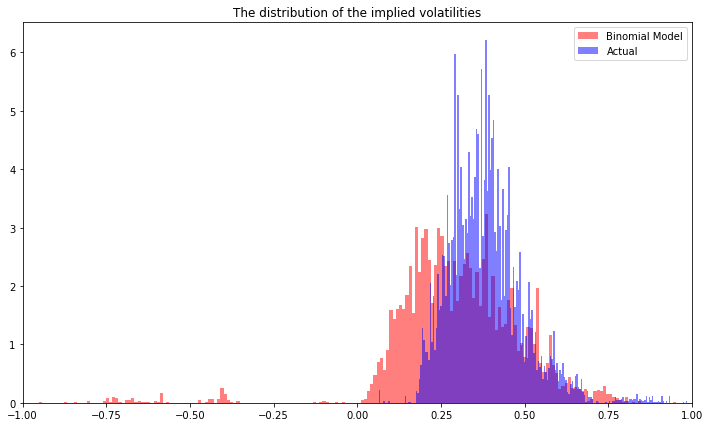}} \\ 
\end{tabular}
\caption{ (left) The scatter plot of implied volatilities calculated from the binomial
model and their corresponding actual values from the dataset (The 2017-2018 US
market); (right) The implied volatilities distributions (The 2017-2018 US market).}
\end{figure}

\subsection{Stock options selected from American markets from 2020-
2021}

\clearpage

\begin{figure}[tph]
\centering
\begin{tabular}{cc}
{\includegraphics[width = 0.5\textwidth, height =
0.55\textwidth]{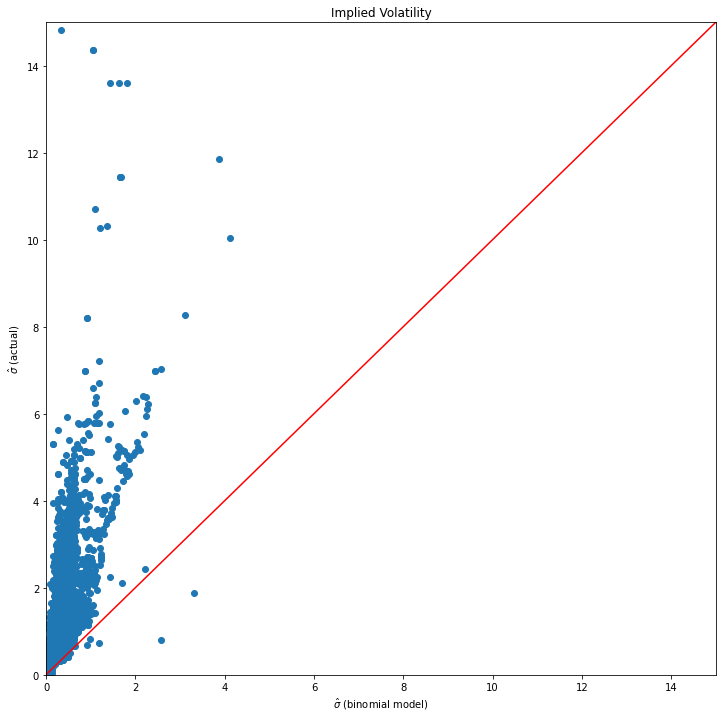}} & {%
\includegraphics[width = 0.5\textwidth,
height = 0.55\textwidth]{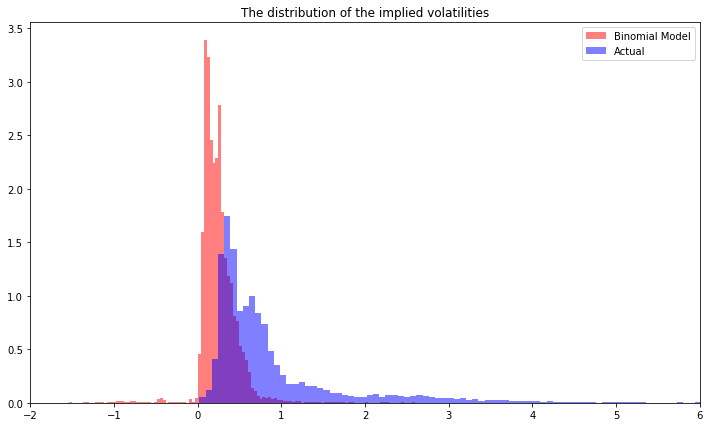}} \\ 
\end{tabular}
\caption{(left) The scatter plot of implied volatilities calculated from the binomial
model and their corresponding actual values from the dataset (The 2020-2021 US
market); (right) The implied volatilities distributions (The 2020-2021 US market).}
\end{figure}

\section{Discussion} In our analysis, we studied the convergence and the rate of convergence of the root-finding methods on both simulated and real data. The results from both datasets showed a consistent interpretation: it takes fewer iterations to compute the implied volatility numerically via automatic differentiation. This result is particularly significant when applied to real data.

Regarding the number of convergences, the secant method performed slightly better than the proposed method. Another interesting observation is that the convergence rate for real-world data, which is around 75\% for both methods, is significantly lower than that for the simulated data, where it exceeds 95\% in both cases. We provide possible explanations for this below:

\begin{enumerate} \item European stock options are typically priced using the Black-Scholes equation rather than the Binomial Model. Although the Binomial Model converges to the Black-Scholes model in continuous time, the implied volatility calculated inversely from different models may not necessarily be the same. This could mean that there is no positive or real root (implied volatility) for the Binomial Model in certain cases. \item The Newton-Raphson method typically finds only one root. Even if the initial value $x_0$ is set close to the average option price, under certain unstable conditions, the method may yield other possible roots, which may not be positive or relevant. This issue can be mitigated by testing the method with different initial values and selecting the most appropriate result. \item Unobservable real-world factors: Unlike the generated dataset, where all variables are controlled, the real stock market is influenced by various uncertainties. \end{enumerate} 

\newpage
\section{Conclusion} In this paper, we present a method to numerically calculate the implied volatility of European stock options. The method consists of two main components: the Binomial Model and a numerical method. The Binomial Model is used to price European stock options. Given data from the exchange market, including the price of the options, we can calculate the implied volatility by inversely solving the model to determine the volatility corresponding to the given option price. The Newton-Raphson method is employed to solve for this volatility. However, the numerical method requires the derivative of the Binomial Model, which can be challenging and tedious to determine analytically or numerically. We propose using Automatic Differentiation, which calculates the derivative of a function immediately after the function's parameters are provided. Finally, we test the proposed method on three datasets: generated data, and real-market data from 2020-2022. The results show that the method performs well with the generated data but its efficiency decreases when tested with real-market data. Additionally, we found that the proposed method converges to the root faster than the secant method.

Practically, the calculated implied volatility can be used to inform option trading decisions. Moreover, this method can be adapted to numerically calculate parameters for other asset pricing models. For example, it can be applied to calculate the implied volatility for American options, which do not follow the Black-Scholes equation but can be evaluated using a modified Binomial Model.

\end{document}